\documentclass[epsf]{article}
\usepackage{amsmath,amsthm,amsfonts,latexsym,amscd}
\usepackage{graphicx}
\title{Self-intersections of Closed Parametrized Minimal Surfaces in Generic Riemannian Manifolds}
\author{John Douglas Moore\\Department of 
Mathematics\\University of California\\
Santa Barbara, CA, USA 93106\\e-mail: moore@math.ucsb.edu}
\date{}

\begin{document}

\maketitle

\begin{abstract}
This article shows that for generic choice of Riemannian metric on a compact manifold $M$ of dimension four, the tangent planes at any self-intersection $p \in M$ of any prime closed parametrized minimal surface in $M$ are not simultaneously complex for any orthogonal complex structure on $M$ at $p$.  This implies via geometric measure theory that $H_2(M;{\mathbb Z})$ is generated by homology classes that are represented by oriented imbedded minimal surfaces. 
\end{abstract}

\section{Introduction}
\label{S:introduction}

The Bumpy Metric Theorem of \cite{M2} (which is also presented as Theorem~5.1.1 in \cite{M}) states that prime parametrized minimal surfaces in a manifold $M$ of dimension at least three are free of branch points when $M$ is given a generic Riemannian metric, and lie on nondegenerate critical submanifolds of the relevant function space.  Once established, it is used to prove the Transversal Crossing Theorem 5.1.2 of \cite{M}, which states that the self-intersections of such minimal surfaces are transverse.  Our goal here is to prove a stronger transversality result when $M$ has dimension four, that the tangent planes are in general position in the sense that they are not simultaneously complex for any orthogonal complex structure on the ambient manifold $M$, which implies via geometric measure theory (GMT) that if $M$ is a compact oriented four-manifold, $H_2(M;{\mathbb Z})$ has a generating set represented by imbedded minimal surfaces.

We discovered this fact soon after publishing the Bumpy Metric Theorem, but deferred publishing it because we believed there should be a better proof which avoids GMT.  The techniques that seem to be needed for this proof have interesting applications that we hope to present elsewhere.  We refer to a recent article of White \cite{W19} for an alternate approach to related transversality results.

If $\hbox{Map}(\Sigma ,M)$ is the space of smooth maps from a surface $\Sigma $ of genus $g$ to $M$ and ${\mathcal T}$ is the Teichm\"uller space of marked conformal structures on compact connected surfaces of genus $g$, a parametrized minimal surface $f : \Sigma \rightarrow M$ can be regarded as a critical point for the {\em Dirichlet energy\/}
\begin{equation} E : \hbox{Map}(\Sigma ,M) \times {\mathcal T} \rightarrow {\mathbb R}, \quad \hbox{defined by} \quad E(f,\omega ) = \frac{1}{2} \int _\Sigma |df|^2dA. \label{E:definitionofenergy}\end{equation}
In this formula, $|df|$ and $dA$ are calculated with respect to some Riemannian metric on $\Sigma $ which lies within the conformal class $\omega \in {\mathcal T}$.  A parametrized minimal surface $f: \Sigma \rightarrow M$ is {\em prime\/} if it is nonconstant and is not a nontrivial cover (possibly branched) of another parametrized minimal surface $f_0: \Sigma _0\rightarrow M$ of lower energy.  By a {\em generic choice of Riemannian metric\/} on $M$ we mean a metric belonging to a countable intersection of open dense subsets of the spaces of $L^2_k$ Riemannian metrics on $M$, as $k$ ranges over the positive integers.

The Bumpy Metric Theorem states that for generic choice of Riemannian metric on a manifold $M$ of dimension at least three, all prime compact oriented parametrized minimal surfaces $f: \Sigma \rightarrow M$ are free of branch points and are as nondegenerate (in the sense of Morse theory) as allowed by the symmetry imposed by the identity component $G$ of the group of conformal automorphisms of $\Sigma $.  If $G$ is trivial, they are Morse nondegenerate in the usual sense, while if $G$ has positive dimension, they lie on nondegenerate critical submanifolds which have the same dimension as $G$.  By a {\em nondegenerate critical submanifold\/} for $F : {\mathcal M} \rightarrow {\mathbb R}$, where ${\mathcal M}$ is a Banach manifold, we mean a submanfold $S \subset {\mathcal M}$ consisting entirely of critical points for $F$ such that the tangent space to $S$ at a given critical point is the space of Jacobi fields for $F$ along $S$.

The Bumpy Metric theorem is a first step towards a partial Morse theory for closed parametrized minimal surfaces in compact four-manifolds with generic Riemannian metrics.  Indeed, the goal of the book \cite{M} was to present foundations for such a partial Morse theory obtained from the known Morse theory for a perturbed $\alpha $-energy of Sacks and Uhlenbeck as it approaches the Dirichlet energy in the limit.  This theory would need to analyze the topological changes that occur in the limit due to bubbling or degeneration of conformal structure.

Example~5.8.2 from \cite{M} shows that even when the metric on the ambient manifold $M$ is generic, branched covers of prime minimal surfaces with nontrivial branch locus need not lie on nondegenerate critical submanifolds.  Thus we need to regard branched covers, together with bubbling and degeneration, as sources of noncompactness (see \S 4.9.3 of \cite{M}) that prevent full Morse inequalities for the Dirichlet energy when the Riemannian metric on the ambient four-manifold is generic.  Yet the topology of the domain of the Dirichlet energy is sufficiently rich that even partial Morse inequalities promise to have striking implications.

For studying nonorientable minimal surfaces, it is convenient to allow the domain of a parametrized surface $f: \Sigma \rightarrow M$ to be nonorientable, and we can study such surfaces by means of their oriented double covers.  For the nonoriented theory one alters the definition of prime to allow nonorientable surfaces.  But Theorem~1 of \cite{M2.5} shows that orientable double covers of nonorientable prime minimal surfaces are also free of branch points and lie on nondegenerate critical submanifolds when the metric on $M$ is generic, so the proof of the Bumpy Metric Theorem for the nonoriented theory (as presented in \cite{M}) implies a corresponding statement for the oriented theory.

We consider the subset
$$\Sigma ^{(s)} = \{ (p_1, \ldots p_s) \in \Sigma ^s : \hbox{$p_i \neq p_j$ when $i\neq j$} \}$$
of the $s$-fold cartesian product $\Sigma ^s$ for $s$ a positive integer, as well as the multidiagonal in the $s$-fold cartesian product $M^s$,
$$\Delta _s = \{ (q_1, \ldots q_s) \in M^s : q_1 = q_2 = \cdots = q_s \}.$$
In accordance with \cite{GG}, Chapter~III, \S 3, we then say that an immersion $f : \Sigma \rightarrow M$ has {\em transversal crossings\/} if for every $s > 1$, the restriction of
$$f^s = f \times \cdots \times f : \Sigma ^s \longrightarrow M^s$$
to $\Sigma ^{(s)}$ is transversal to $\Delta _s$.  Thus if $\Sigma $ is a compact surface and $M$ has dimension at least five, an immersion with transversal crossings is a one-to-one immersion and hence an imbedding, while if $M$ has dimension four, such an immersion has only double points and the intersections at double points are transverse.

We can regard the Transversal Crossing Theorem 5.1.2 of \cite{M} as stating that having transversal crossings is a {\em generic property\/} of parametrized minimal surfaces in four-manifolds, a generic property being one that holds when the Riemannian metric on the ambient manifold $M$ is generic.  For the development of our theory, we need a procedure for establishing additional generic properties which might be useful in establishing our desired partial Morse inequalities.  Fortunately, the Transversal Density Theorem~19.1 of Abraham and Robbin \cite{AR} provides a relatively straightforward procedure for checking that properties are generic.  As an application of the Abraham-Robbin theory, we will prove:

\vskip .1 in
\noindent
{\bf Theorem~1.} {\sl Suppose that $M$ is a compact manifold of dimension four.  Then for a generic choice of Riemannian metric on $M$, at any self-intersection point of a transversely immersed minimal surface $f : \Sigma \rightarrow M$, the tangent planes are in general position with respect to the metric, that is, they are not simultaneously complex for any orthogonal complex structure on the tangent space at a point of self-intersection.}

\vskip .1 in
\noindent
This theorem is likely needed as a step towards establishing the partial Morse inequalities for minimal surfaces in four-manifolds that we believe exist.

Recall that an {\em orthogonal complex structure\/} on the tangent space $T_pM$ to a Riemannian manifold $M$ at a point $p \in M$ is an endomorphism
\begin{multline*} J : T_pM \rightarrow T_pM \quad \hbox{such that} \quad J^2 = -I \\ \hbox{and} \quad \langle Jv, Jw \rangle = \langle v,w \rangle , \quad \hbox{for $v,w \in T_pM$,}\end{multline*}
where $\langle \cdot , \cdot \rangle $ denotes the Riemannian metric on $M$.  When $M$ is four-dimensional, we say that $J$ is positively or negatively oriented depending on whether the corresponding {\em K\"ahler form\/}, defined by
$$\theta (v,w) = \langle Jv, w \rangle, \quad \hbox{for} \quad v,w \in T_pM,$$
is self-dual or anti-self-dual.  A globally defined orthogonal complex structure on the Riemannian manifold $M$ is a family of orthogonal complex structure $p \mapsto J(p)$ which depends smoothly on the parameter $p$ which ranges over $M$.

The generic condition on the tangent planes enables us to use a result of Frank Morgan \cite{Mor} to show that if $f: \Sigma \rightarrow M$ is a surface of genus $g$ which minimizes area in some homology class, and $f$ has points of self-intersection, then one of the self-intersections can be removed by surgery, producing a surface of larger genus and smaller area in the same homology class.

\vskip .1 in
\noindent
{\bf Theorem~2.} {\sl Suppose that $M$ is a compact oriented manifold of dimension at least four with a generic choice of Riemannian metric.  Then each nonzero element of $H_2(M;{\mathbb Z})$ is represented by a collection of disjoint connected minimal surfaces, each of which is either imbedded or a nontrivial branched cover of an imbedded minimal surface.}

\vskip .1 in
\noindent
Assuming Theorem~1, we can prove Theorem~2 as follows.  Results of Almgren and Chang \cite{Chang} (see the Main Regularity Result on page~72 of \cite{Chang}) imply that any homology class is represented by an area minimizing integral current which arises from a smooth submanifold except for possible branch points and self-intersections.  This can be represented by a finite collection of parametrized minimal surfaces, each of which is either prime or a branched cover of a prime minimal surface.  Let $ f_i : \Sigma _i \rightarrow M $ for $1 \leq i \leq k$ be the underlying prime minimal surfaces, where each $\Sigma _i$ is connected.  When the metric is generic, it follows from the Bumpy Metric Theorem that each such $f_i$ is free of branch points, while when the dimension of $M$ is at least five, it follows from the Transversal Crossing Theorem that there are no self-intersections, or intersections between different components.  When the dimension of $M$ is four, the Transversal Crossing Theorem and Theorem~1 imply that intersections are transverse, and that at any self-intersection the two tangent planes of $f_i$ cannot be simultaneously complex for any orthogonal complex structure at the point of intersection.  It therefore follows from Theorem~2 of \cite{Mor} that if any $f_i$ has nontrivial self-intersections, one of the self-intersections can be eliminated with a decrease in area, thereby contradicting the fact that the current is area minimizing.  Thus the $f_i$'s must be imbeddings.  Similarly, the area could be decreased if the images of different $f_i$'s were not mutually disjoint, again contradicting area minimization.  QED

\vskip .1 in
\noindent
Theorem~2 is related to an earlier result of Brian White \cite{W} which treats unoriented surfaces.

\section{Proof of Theorem~1}
\label{S:proof1}

In preparation, we recall some well-known facts about manifolds of maps.  From the survey of Eells \cite{Ee}, we find an elegant combination of the $\alpha$- and $\omega $-lemmas of Abraham and Smale (as presented in \cite{Ab}) which asserts that
$$\Phi : C^{k+r}(M,N) \times C^k(S,M) \rightarrow C^k(S,N), \qquad \Phi (f,g) = f \circ g, \quad r, k \in {\mathbb N} \cup \{ 0 \},$$
is a $C^r$ map when $S$, $M$ and $N$ are compact smooth manifolds.  A direct consequence (in the case where $S$ is a point) is the fact that the evaluation map
$$\hbox{ev} : C^{r}(M ,N) \times M \rightarrow N, \quad \hbox{ev}(f,p) = f(p)$$
is $C^r$, a fact used in the Abraham-Robbin approach to transversality theory, as explained in \cite{AR}.  If ${\mathcal M}$ is a Banach manifold which satisfies the second axiom of countability, then a map $\rho : {\mathcal M} \rightarrow C^{r}(M ,N)$ is said to be a {\em $C^r$-representation\/} if
$$\hbox{ev}_\rho : {\mathcal M} \times M \rightarrow N, \quad \hbox{defined by} \quad \hbox{ev}_\rho(f,p) = \rho(f)(p)$$
is $C^r$.  For example, the Sobolev Lemma shows that the inclusion
$$\rho : L^2_k(\Sigma ,M) \rightarrow C^r(\Sigma ,M)$$
is a $C^{k-2}$-representation when $\Sigma $ is compact of dimension two and $k \geq 2$.  Here $L^2_k(\Sigma ,M)$ is the completion of the space of $C^\infty $ maps from $\Sigma $ to $M$ with respect to the Sobolev $L^2_k$ norm. 

For simplicity, let us assume first that $\Sigma $ is connected and has genus at least two.  We need to be precise about the regularity of the maps in question, so we will replace $\hbox{Map}(\Sigma ,M)$ and $\hbox{Met}(M)$ by
$$L^2_k(\Sigma ,M) \quad \hbox{and} \quad \hbox{Met}^2_{\ell}(M),$$
where $\hbox{Met}^2_{\ell}(M)$ is the $L^2_{\ell }$ completion of the space of $C^\infty $ Riemannian metrics on $M$, with $\ell \geq k-1$.  It is shown in \S 5.4 of \cite{M} that
\begin{multline}{\mathcal P}^{\emptyset} = \{ (f,\omega ,g) \in L^2_k(\Sigma ,M) \times {\mathcal T} \times \hbox{Met}^2_{\ell}(M) :  \hbox{ $(f,\omega)$ is a prime} \\ \hbox{immersed conformal harmonic map for the metric $g$ }  \}. \label{E:defofP} \end{multline}
is a smooth submanifold, and that the projection
$$\pi_2 : {\mathcal P}^{\emptyset} \longrightarrow \hbox{Met}(M)$$
is a Fredholm map of Fredholm index zero.  (The superscript on ${\mathcal P}^{\emptyset}$ indicates that its elements have trivial branch locus.)  Let
$${\mathcal P}^{\emptyset}_0 = \{ (f,\omega ,g) \in  {\mathcal P}^{\emptyset} : \hbox{ $f$ is Morse nondegenerate } \}.$$
When $k$ and $\ell > k$ are large, regularity results show that elements of ${\mathcal P}^{\emptyset}$, being critical points of the energy, already have higher regularity than assumed:
$${\mathcal P}^{\emptyset} \subseteq L^2_\ell (\Sigma ,M) \times {\mathcal T} \times \hbox{Met}^2_{\ell}(M).$$
Thus we can assume that the evaluation map has as many derivatives as we want by choosing $\ell $ large.  If $(f,\omega ,g)$ is an element of ${\mathcal P}^{\emptyset}_0$, it follows from the inverse function theorem, that there is an open neighborhood ${\mathcal U}$ of $(f,\omega ,g)$ in ${\mathcal P}^{\emptyset}_0$ an open neighborhood ${\mathcal V}$ of $g$ in $\hbox{Met}(M)$ and a diffeomorphism
\begin{equation} \sigma : {\mathcal V} \rightarrow {\mathcal U} \quad \hbox{such that} \quad \pi _2\circ \sigma = \hbox{id}_{\mathcal V}. \label{E:sigma}\end{equation}
If $\pi _0 : {\mathcal P}^{\emptyset}_0 \rightarrow L^2_\ell (\Sigma ,M)$ is the projection on the first factor, we can regard
\begin{equation} \rho = \pi _2 \circ \sigma : {\mathcal V} \rightarrow L^2_\ell (\Sigma ,M) \subseteq C^{\ell - 2}(\Sigma ,M) \label{E:rho}\end{equation}
as a $C^{\ell - 2}$-representation.  Since the manifolds in question have countable base and the projection $\pi _0$ is locally proper by Theorem 1.6 of \cite{Sm2}, we can choose a sequence
$$\{ \sigma _i : {\mathcal V}_i \rightarrow {\mathcal P}^{\emptyset}_0 : i \in {\mathbb N} \} \quad \hbox{such that} \quad \bigcup \sigma _i\left({\mathcal V}_i\right) = {\mathcal P}^{\emptyset}_0.$$
Thus it suffices to show that for any section $\sigma $ constructed as above and for a residual subset of $g \in {\mathcal V}$, $\rho (g)$ is an immersion such that at any point of self-intersection the tangent planes to $\rho (g)$ at the intersection point are not simultaneously complex for any choice of orthogonal complex structure on the tangent space to $M$.

To prove this, it is convenient to construct Gauss lifts of (\ref{E:rho}) to the six-dimensional twistor spaces, two-sphere bundles over $M$, which were associated to $M$ by Eells and Salamon \cite{ES}.  If $M$ is spin these twistor spaces are defined in terms of the spin bundles $W_-$ and $W_+$ over $M$.  These are $SU(2)$ bundles over $M$ which give a factorization of the complexified tangent bundle,
$$TM \otimes {\mathbb C} \cong \hbox{Hom}(W_+,W_-) \cong W_- \otimes W_+^*,$$
as described in \S 2.2 of \cite{M0}.  The Eells-Salamon twistor spaces are then defined to be the projective bundles ${\mathbb P}(W_+)$ and ${\mathbb P}(W_-)$ over $M$.  If $M$ is not spin, we can still define $W_-$ and $W_+$ as \lq\lq virtual" $SU(2)$ bundles over $M$ (see \S 2.3 of \cite{M0}), and this still allows us to construct the twistor spaces.  Given a parametrized minimal surface $f : \Sigma \rightarrow M$, Eells and Salamon construct two corresponding {\em Gauss lifts\/}
$$f_+ : \Sigma \rightarrow {\mathbb P}(W_+) \quad \hbox{and} \quad f_- : \Sigma \rightarrow {\mathbb P}(W_-),$$
and show that these lifts are pseudoholomorphic with respect to suitable orthogonal complex structures globally defined on ${\mathbb P}(W_+)$ and ${\mathbb P}(W_-)$.  We can think of $f_+$ and $f_-$ as defining families of orthogonal complex structures at points of the image of $f$.  To prove Theorem~1 we need to show that when the metric on $M$ is generic, the orthogonal complex structures thus assigned to distinct points $p$ and $q$ of $\Sigma $ disagree whenever $f(p) = f(q)$, for every closed parametrized minimal surface $f : \Sigma \rightarrow M$.

In the notation we used earlier, we construct Gauss lifts
$$\rho _+ : {\mathcal V} \rightarrow L^2_{\ell -1}(\Sigma ,{\mathbb P}(W_+)), \quad \rho _- : {\mathcal V} \rightarrow L^2_{\ell -1}(\Sigma ,{\mathbb P}(W_-))$$
of the map $\rho $ of (\ref{E:rho}), where now the projective spinor bundles ${\mathbb P}(W_+)$ and ${\mathbb P}(W_-)$ vary smoothly with the metric.  These are $C^{\ell - 3}$-representations of ${\mathcal V}$.  It then suffices to show that for generic choice of $g \in{\mathcal V}$, the Gauss lifts $\rho _+(g)$ and $\rho _-(g)$ of $\rho (g)$ are imbeddings into ${\mathbb P}(W_+)$ and ${\mathbb P}(W_-)$ respectively.  By taking the Cartesian product of $\rho _+(g)$ and $\rho _-(g)$ with themselves, we define $C^{\ell - 3}$-representations
\begin{multline} \rho _{+,2} : {\mathcal V} \rightarrow L^2_{\ell -1}(\Sigma \times \Sigma - \Delta ,{\mathbb P}(W_+) \times {\mathbb P}(W_+)), \\ \rho _{-,2} : {\mathcal V} \rightarrow L^2_{\ell -1}(\Sigma \times \Sigma - \Delta ,{\mathbb P}(W_-) \times {\mathbb P}(W_-)),\label{E;rhoplusminus2}\end{multline}
where $\Delta $ is the diagonal.  This puts us in a position to apply the Transversal Density Theorem~19.1 of \cite{AR}.

To prove our theorem, it suffices to show that if $\Delta _+$ and $\Delta _-$ are the diagonals in ${\mathbb P}(W_+) \times {\mathbb P}(W_+)$ and ${\mathbb P}(W_-) \times {\mathbb P}(W_-)$, respectively, then
\begin{multline*} {\mathcal V}_{\Delta _+} = \{ g \in {\mathcal V} : \hbox{ $\rho _+(g)$ is transversal to $\Delta _+$ } \}, \\  {\mathcal V}_{\Delta _-} = \{ g \in {\mathcal V} : \hbox{ $\rho _-(g)$ is transversal to $\Delta _-$ } \}\end{multline*}
are residual subsets of ${\mathcal V}$, and the Transversal Density Theorem ensures that this will be the case if the maps
\begin{multline*}\hbox{ev} \circ \rho _+: {\mathcal V} \times (\Sigma \times \Sigma - \Delta ) \rightarrow {\mathbb P}(W_+) \times {\mathbb P}(W_+) \\ \hbox{and} \quad \hbox{ev} \circ \rho _-: {\mathcal V} \times (\Sigma \times \Sigma - \Delta ) \rightarrow {\mathbb P}(W_-) \times {\mathbb P}(W_-)\end{multline*}
are transversal to $\Delta _+$ and $\Delta _-$ respectively.


In other words, assuming that $M$ has dimension four, we need to construct a variation of the metric which puts a given intersection into general position, the two intersecting planes not being simultaneously complex for any orthogonal complex structure.  

Suppose therefore that $p$ and $q$ are distinct points of $\Sigma $ with $f(p) = f(q)$, and let $V_1$ and $V_2$ be disjoint open neighborhoods of $p$ and $q$ within $\Sigma $.  We construct coordinates $(u^1, u^2, u^3, u^4)$ on a neighborhood $U$ of $f(p)$ in $M$ so that
\begin{enumerate}
\item $u^i(f(p)) = 0$,
\item $f(V_1) \cap U$ is described  by the equations $u^3 = u^4 = 0$,
\item $f(V_2) \cap U$ is described  by the equations $u^1 = u^2 = 0$,
\item $f^*\langle \cdot, \cdot \rangle|V_1 = \lambda _1^2 ((dx^1)^2 + (dx^2)^2)$, where $x^a = u^a \circ f$, and
\item $f^*\langle \cdot, \cdot \rangle|V_2 = \lambda _2^2 ((dx^3)^2 + (dx^4)^2)$, where $x^r = u^r \circ f$. 
\end{enumerate}
Here our index conventions will restrict indices to the limits
$$1 \leq a,b \leq 2, \qquad 3 \leq r,s \leq 4, \qquad 1 \leq i,j,k \leq 4.$$
Let $g_{ij}$ be the components of the metric in these coordinates, so that
$$g_{ab} = \lambda _1^2 \delta _{ab} \quad \hbox{along $f(V_1)$} \qquad \hbox{and} \qquad g_{rs} = \lambda _2^2 \delta _{rs} \quad \hbox{along $f(V_2)$,}$$
where $\delta _{ab}$ and $\delta _{rs}$ denote Kronecker deltas.  We assume that at the intersection point, $f_*(T_p\Sigma )$ and $f_*(T_q\Sigma )$ are simultaneously complex for some orthogonal complex structure on $TM$.  (After reordering $u^3$ and $u^4$ if necessary, we can then assume without loss of generality that $g_{13} = g_{24}$ and $g_{14} = - g_{23}$.)

If we define the Christoffel symbols in terms of the metric
$$\Gamma _{k,ij} = \frac{1}{2} \left( \frac{\partial g_{ki}}{\partial u^i} + \frac{\partial g_{kj}}{\partial u^i} - \frac{\partial g_{ij}}{\partial u^k}\right), \quad \Gamma ^k_{ij} = \sum g^{kl} \Gamma _{l,ij},$$
the fact that $f$ is harmonic is expressed by the equations
\begin{equation} \Gamma ^k_{11} + \Gamma ^k_{22} = 0 \quad \hbox{along $f(V_1)$,} \qquad \Gamma ^k_{33} + \Gamma ^k_{44} = 0 \quad \hbox{along $f(V_2)$.} \label{E:eqforharmonic}\end{equation}
 
We will construct a variation in the metric $(\dot g_{ij})$ such that $\dot g_{ab} = 0 = \dot g_{rs}$ and the equations (\ref{E:eqforharmonic}) continue to hold.  The resulting variation $\dot \Gamma _{k,ij}$ in the Christoffel symbols will then satisfy the equations
$$\dot \Gamma _{b,aa} = 0, \quad \dot \Gamma _{r,aa} = \frac{\partial \dot g_{ra}}{\partial u^a}, \quad \dot \Gamma _{s,rr} = 0, \quad \dot \Gamma _{a,rr} = \frac{\partial \dot g_{ra}}{\partial u^r}.$$
Thus we want to arrange that
\begin{equation}\sum _a\frac{\partial \dot g_{ra}}{\partial u^a} = 0 \quad \hbox{along $f(V_1)$,} \quad \hbox{and} \quad \sum _r\frac{\partial \dot g_{ra}}{\partial u^r} = 0 \quad \hbox{along $f(V_2)$.}\label{E:necesaary}\end{equation}
If we construct a smooth function $h : U \rightarrow {\mathbb R}$ and then set
$$\begin{pmatrix} \dot g_{13} & \dot g_{14} \cr \dot g_{23} & \dot g_{24} \end{pmatrix} = \begin{pmatrix} \frac{\partial ^2h}{\partial u^2\partial u^4} & - \frac{\partial ^2h}{\partial u^2\partial u^3} \cr - \frac{\partial ^2h}{\partial u^1\partial u^4} & \frac{\partial ^2h}{\partial u^1\partial u^3}\end{pmatrix},$$
we find that the equations (\ref{E:necesaary}) are satisfied.  We can choose such a function which has compact support within $U$, and for which
$$\begin{pmatrix} \dot g_{13} & \dot g_{14} \cr \dot g_{23} & \dot g_{24} \end{pmatrix} (f(p))$$
is arbitrary.  The resulting metric perturbation will preserve conformality and minimality of $f$ as required, yet can be chosen so that after perturbation $f_*(T_p\Sigma )$ and $f_*(T_q\Sigma )$ will not be simultaneously complex for any orthogonal complex structure on $T_{f(p)}M = T_{f(q)}M$.  This finishes the proof of Theorem~1 in the case where $\Sigma $ is connected and has genus at least two.

The additional complication we need to handle when the genus is zero or one is that the group $G$ of complex automorphisms of $\Sigma $ which are homotopic to the identity is nontrivial.  One approach to treating this case is to restrict to submanifolds of $L^2_k(\Sigma , M)$ which have codimension equal to the dimension of $G$, yet represent all immersions up to the action of $G$.

Suppose first that $\Sigma$ is the two-sphere $S^2$, in which case $G = PSL(2,{\mathbb C})$.  If $N$ is a compact codimension two submanifold of $M$ with boundary $\partial N$, we let
\begin{multline}{\mathcal U}(N) = \{ f \in L^2_k(S^2, M) : \hbox{ $f$ does not intersect $\partial N$ and has nonempty }\\ \hbox{ transversal intersection with the interior of $N$ } \}, \label{E:un1} \end{multline}
which is an open subset of $L^2_k(S^2, M)$.  Given three disjoint compact codimension two submanifolds with boundary, say $Q$, $R$ and $S$, we let
\begin{equation}{\mathcal U}(Q,R,S) = {\mathcal U}(Q) \cap {\mathcal U}(R) \cap {\mathcal U}(S), \label{E:un2} \end{equation}
also an open subset of $L^2_k(S^2, M)$.  We cover $L^2_k(S^2, M)$ with a countable collection of sets ${\mathcal U}(Q_j,R_j,S_j)$ defined by a sequence $j \mapsto (Q_j,R_j,S_j)$ of triples of such codimension two submanifolds, where $j \in {\mathbb N}$.  We choose three points $q$, $r$ and $s$ in $S^2$  and let
\begin{equation}{\mathcal F}_j(S^2 ,M) = \{ f \in {\mathcal U}(Q_j,R_j,S_j) : f(q) \in Q_j, f(r) \in R_j, f(s) \in S_j  \},\label{E:defoffiS2}\end{equation}
noting that ${\mathcal F}_j(S^2 ,M)$ meets each $PSL(2,{\mathbb C})$-orbit in ${\mathcal U}(Q_j,R_j,S_j)$ in a finite number of points.  It follows from the Sobolev imbedding theorem and smoothness of the evaluation map on the space of $C^{k-2}$ maps (for example by Proposition~2.4.17 of \cite{AMR}) that the evaluation map
$$\hbox{ev} : C^k(S^2, M) \times S^2 \longrightarrow M, \qquad \hbox{defined by} \qquad \hbox{ev}(f,p) = f(p),$$
is $C^{k-2}$.  Thus when $k$ is large, we can regard ${\mathcal F}_j(S^2,M)$ as a submanifold of $\hbox{Map}(S^2, M)$ of codimension six with tangent space
\begin{multline*}T_f{\mathcal F}_j(S^2,M) = \{ \hbox{sections $X$ of $f^*TM$} : \\ X(q) \in T_{f(q)}Q_j, X(r) \in T_{f(r)}R_j, X(s) \in T_{f(s)}S_j \}.\end{multline*}

In the case where $\Sigma $ is the torus and $G = S^1 \times S^1$, we need fix only one point to break the symmetry.  We define ${\mathcal U}(N)$ by (\ref{E:un1}) with $S^2$ replaced by $T^2$ and choose a sequence $j \mapsto N_j$ of smooth compact codimension two submanifolds of $M$ such that ${\mathcal U}(N_j)$ cover $\hbox{Map}(T^2, M)$.  We choose a base point $q \in T^2$ and let
\begin{equation} {\mathcal F}_j(T^2 ,M)= \{ f \in {\mathcal U}(N_j) : f(q) \in N_j  \}. \label{E:defoffiT2}\end{equation}

We now simply carry through the earlier argument with $L^2_k(\Sigma , M)$ with a collection of spaces ${\mathcal F}_j(\Sigma ,M)$ which cover $L^2_k(\Sigma , M)$.

Finally, we need to consider the case in which $\Sigma $ has several components, and in particular we need to consider intersections between two different components.  But the modifications needed to treat this case are straightforward and we can safely leave them to the reader.  QED

\vskip .1in
\noindent
{\bf Acknowledgement}.  The author thanks the referees for several helpful suggestions.

\vskip .1in
\noindent
{\bf Data availability statement}.  Data sharing is not applicable to this article because no data sets were used.

\end{document}